\documentclass[12pt]{article}

\usepackage{amsmath, epsfig, cite}
\usepackage{amssymb}
\usepackage{amsfonts}
\usepackage{latexsym}
\usepackage{graphicx}

\newtheorem{thm}{Theorem}[section]

\newtheorem{cor}[thm]{Corollary}

\newtheorem{rem}[thm]{\it Remark}

\newcommand{\qed}{{\hfill\rule{4pt}{7pt}}}

\numberwithin{equation}{section}

\makeatletter \@addtoreset{equation}{section} \makeatother

\begin{document}
\rule{0cm}{3cm}
\begin{center}
{\Large\bf On the Bilateral Series $_2\psi_2$}
\end{center}
\vskip 2mm \centerline{Vincent Y. B. Chen$^1$, William Y. C.
Chen$^2$ and Nancy S. S. Gu$^3$ }

\begin{center}
Center for Combinatorics, LPMC\\
Nankai University, Tianjin 300071\\
People's Republic of China\\

\vskip 2mm
 Email:  $^1$ybchen@mail.nankai.edu.cn, $^2$chen@nankai.edu.cn,
 $^3$gu@nankai.edu.cn
\end{center}

\begin{center}
{\bf Abstract}
\end{center}

{\small We obtain a formula which reduces the evaluation of a
$_2\psi_2$ series to two $_2\phi_1$ series. In some sense, this
identity may be considered as a companion  of Slater's formulas. We
also find that a two-term ${}_2\psi_2$ summation formula due to
Slater can be  derived from a unilateral summation formula of
Andrews by bilateral extension and parameter augmentation. }

\vskip 3mm



\section{Introduction}

It is well known that many   bilateral basic hypergeometric
identities can be derived from unilateral identities. Using Cauchy's
method \cite{Bailey1, JouhetandSchosser, Schlosser, Schlosser1} one
may obtain bilateral basic hypergeometric identities from
terminating unilateral identities. Starting with nonterminating
unilateral basic hypergeometric series, Chen and Fu  \cite{chenfu}
developed a method to construct semi-finite forms by shifting  the
summation index by $m$. Then the bilateral summations are
consequences of the semi-finite forms  by  letting $m$ tend to
infinity. We call this method bilateral extension.  In this paper we
use bilateral extensions of a ${}_3\phi_2$ series and an identity of
Andrews \cite{Andrews} to study the bilateral series ${}_2\psi_2$:
\begin{equation}
{}_{2}\psi_2\left[\begin{array}{cc}
a,&b\\
c,& d
\end{array};q,z\right].\label{2psi2-gene}
\end{equation}

The above ${}_2\psi_2$ series is closely related to
 the question of finding a
$q$-extension of Dougall's bilateral hypergeometric series summation
formula \cite{Dougall}:
\begin{eqnarray}\label{Dougall}
\sum_{k=-\infty}^{\infty}\frac{(a)_k(b)_k}{(c)_k(d)_k}=
\frac{\Gamma(c)\Gamma(d)\Gamma(1-a)\Gamma(1-b)\Gamma(c+d-a-b-1)}
{\Gamma(c-a)\Gamma(c-b)\Gamma(d-a)\Gamma(d-b)},
\end{eqnarray}
where $\text{Re}(c+d-a-b-1)>0$, $(a)_k=a(a+1)\cdots(a+k-1)$, $k=1,
2, \cdots$, $(a)_0=1$ and $(a)_k=(-1)^k/(1-a)_{-k}$ when $k$ is a
negative integer.

Bailey \cite{Bailey2} first suggested that there did not exist any
$q$-extension of \eqref{Dougall}. Since (\ref{Dougall}) is an
extension of the Gauss ${}_2F_1$ summation formula, one naturally
expects that a $q$-analogue of (\ref{Dougall}) should be concerned
with the following series:
\begin{equation}
{}_{2}\psi_2\left[\begin{array}{cc}
a,&b\\
c,& d
\end{array};q,{cd\over abq} \right].\label{2psi2gauss}
\end{equation}
Clearly, when $c$ or $d$ equals $q$, \eqref{2psi2gauss} reduces to
the $q$-Gauss sum
 \cite[Appendix
II.8]{Gasper-Rahman}:
\begin{eqnarray}
{}_{2}\phi_{1}\left[\begin{array}{ccc}
a,&b\\
&c
\end{array};q,\frac{c}{ab}\right]=
\frac{(c/a,c/b;q)_\infty}
{(c,c/ab;q)_\infty},\quad\quad\quad|c/ab|<1.\label{Gauss}
\end{eqnarray}
Even for the above series (\ref{2psi2gauss}), Gasper  \cite{Gasper}
pointed out that one could not use analytic continuation to derive an infinite product representation.

On the other hand, many  results on
the bilateral $_2\psi_2$ series \eqref{2psi2-gene} have been obtained. In
\cite{Bailey2}, Bailey found several transformation formulas for the
$_2\psi_2$ series \eqref{2psi2-gene}. Later, Slater obtained a
general transformation formula for an $_r\psi_r$ series in \cite{slater0} based
on Sears'  transformation on the $_{r+s+1}\phi_{r+s}$
series in \cite{Sears}  subject to suitable substitutions and the following  relation
\begin{equation}
\sum_{n=0}^{\infty}f(n)=\sum_{n=-\infty}^{-1}f(-n-1)
\end{equation}
to combine two unilateral series to form a bilateral series.

Gasper and Rahman \cite{Gasper-Rahman} have shown that based on
Slater's transformation formula, one could obtain two expansions of
an $_r\psi_r$ series in terms of $r$ $_r\phi_{r-1}$ series \cite[Eq.
(5.4.4), (5.4.5)]{Gasper-Rahman}. When $r=2$, they become
\begin{eqnarray} \label{22-1}
&&{}_2\psi_2\left[\begin{array}{cc}
a, &b\\
c, & d
\end{array};q,z\right]=\frac{a(q,qa/b,c/a,d/a,az,q/az,qb,1/b;q)_{\infty}}
{(a/b,qb/a,c,d,q/a,q/b,z,q/z;q)_{\infty}}\nonumber \\
&&\quad\quad\quad\quad\times{}_2\phi_1\left[\begin{array}{cc}
qa/c, &qa/d\\
&qa/b
\end{array};q,\frac{cd}{abz}\right]+\text{idem}(a;b)
\end{eqnarray}
and \begin{eqnarray} \label{22-2}
&&{}_2\psi_2\left[\begin{array}{cc}
a, &b\\
c, & d
\end{array};q,z\right]=\frac{q}{c}\frac{(q,c/a,c/b,abz/dq,dq^2/abz,q/d;q)_{\infty}}
{(c,c/d,q/a,q/b,abz/cd,qcd/abz;q)_{\infty}}\nonumber \\
&&\quad\quad\quad\quad\times{}_2\phi_1\left[\begin{array}{cc}
qa/c, &qb/c\\
&qd/c
\end{array};q,z\right]
+\text{idem}(c;d),
\end{eqnarray}
where the symbol ``idem($a;b$)'' after an expression means that the
preceding expression is repeated with $a$ and $b$ interchanged.

Setting $d=q$, \eqref{22-1} reduces to a three-term transformation
formula \cite[Appendix III.32]{Gasper-Rahman} for the ${}_2\phi_1$
series:
\begin{equation}\label{three}
{}_2\phi_1\left[\begin{array}{cc}
a, &b\\
& c
\end{array};q,z\right]=\frac{(b,c/a,az,q/az;q)_{\infty}}{(c,b/a,z,q/z;q)_{\infty}}
{}_2\phi_1\left[\begin{array}{cc}
a, &aq/c\\
& aq/b
\end{array};q,\frac{cq}{abz}\right]+\text{idem}(a;b).
\end{equation}
However, it should be noted that when $c$ or $d$ equals $q$, \eqref{22-2} does not lead to any
nontrivial identity.

The first result of this paper is to give a new formula for the
$_2\psi_2$ series \eqref{2psi2-gene} in terms of two $_2\phi_1$
series which is different from Slater's formulas \eqref{22-1} and
\eqref{22-2}. It reduces to a different three-term transformation
formula \eqref{three-1} when $c=q$ compared with the three-term
transformation formula (\ref{three}) deduced by Slater's
transformation. Moreover, this identity may be considered as a
companion of Slater's formulas \eqref{22-1} and \eqref{22-2}.  Note
that Slater's formulas do not seem to imply the special cases that
can be deduced from our formula except for Ramanujan's ${}_1\psi_1$
summation formula \cite[Appendix II.29]{Gasper-Rahman}. As a
consequence, our formula yields a two-term closed product form for
the ${}_2\psi_2$ series:
\begin{eqnarray}\label{2psi2-2}
&&{}_{2}\psi_{2}\left[\begin{array}{ccc}
b,&c\\
aq/b,&aq/c
\end{array};q,-\frac{aq}{bc}\right]=
\frac{(-b,aq/bc,-q/b,b/a,q;q)_\infty(aq^2/c^2;q^2)_\infty}
{(aq/c,-1,q/c,q/b,-aq/bc;q)_\infty(b^2/a;q^2)_\infty}\nonumber\\
&&\quad\quad\quad\quad\quad\quad\quad\quad
+\frac{(aq/bc,b,-aq/b,-b/a,q;q)_\infty(aq^2/c^2;q^2)_\infty}
{(aq/b,aq/c,-1,-aq/bc,q/c;q)_\infty(b^2/a;q^2)_\infty}.
\end{eqnarray}
For comparison, we recall the known formula for the well-poised
$_2\psi_2$ series \cite[Appendix II.30]{Gasper-Rahman}:
\begin{eqnarray}
{}_{2}\psi_{2}\left[\begin{array}{ccc}
b,&c\\
aq/b,&aq/c
\end{array};q,-\frac{aq}{bc}\right]=
\frac{(aq/bc;q)_\infty(aq^2/b^2,aq^2/c^2,q^2,aq,q/a;q^2)_\infty}
{(aq/b,aq/c,q/b,q/c,-aq/bc;q)_\infty}.\label{2psi21term}
\end{eqnarray}

Let us turn our attention back to Dougall's formula. As pointed out
by Askey \cite{Askey}, Bailey seemed to have been partly right
concerning his opinion towards the $q$-extension of Dougall's
formula. According to Askey \cite{Askey},  in certain sense the
following $q$-extension of Cauchy's beta integral was similar to a
$q$-extension of Dougall's formula:
\begin{equation}\label{integral}
\int_{-\infty}^{\infty}\frac{(ct,-dt;q)_{\infty}}{(at,-bt;q)_{\infty}}\text{d}_qt
=2\frac{(1-q)(c/a,d/b,-c/b,-d/a,ab,q/ab;q)_{\infty}(q^2;q^2)_{\infty}^2}{(cd/abq,q;q)_{\infty}
(a^2,q^2/a^2,b^2,q^2/b^2;q^2)_{\infty}}.
\end{equation}
In fact, this integral can be recast as  a two-term summation
formula for the $_2\psi_2$ series \eqref{2psi2gauss}:
\begin{eqnarray}
&&\frac{(c,-d;q)_{\infty}}{(a,-b;q)_{\infty}}
{}_{2}\psi_{2}\left[\begin{array}{ccc}
a,&-b\\
c,&-d
\end{array};q,q\right]+\frac{(-c,d;q)_{\infty}}{(-a,b;q)_{\infty}}
{}_{2}\psi_{2}\left[\begin{array}{ccc}
-a,&b\\
-c,&d
\end{array};q,q\right]\nonumber\\[5pt]
&&\quad\quad\quad=2\frac{(1-q)(c/a,d/b,-c/b,-d/a,ab,q/ab;q)_{\infty}(q^2;q^2)_{\infty}^2}{(cd/abq,q;q)_{\infty}
(a^2,q^2/a^2,b^2,q^2/b^2;q^2)_{\infty}}.
\end{eqnarray}

As observed by Ismail and Rahman \cite{Ismail-Rahman}, the above
two-term summation formula is a special of a transformation formula
due to Slater \cite{slater0}. When $r=2$, by substitutions and the
$q$-Gauss sum \eqref{Gauss}, Slater's general transformation on the
$_r\psi_r$ series reduces to the following two-term summation
formula:
\begin{eqnarray}\label{Slater}
&&\frac{(c/ef,qef/c,q,q/a,q/b,c/a,c/b;q)_{\infty}}{(e,f,q/e,q/f,c/ab;q)_{\infty}}
=\frac{q}{e}\frac{(c/qf,q^2f/c,e/a,e/b,qc/e,q^2/e;q)_{\infty}}{(e,q/e,e/f,qf/e;q)_{\infty}}
\nonumber \\[5pt]&&\quad\quad\quad\quad\times{}_2\psi_2\left[\begin{array}{cc}
e/c, &e/q\\
e/a,&e/b
\end{array};q,q\right]+\text{idem}(e;f).
\end{eqnarray}

The second result of this paper is concerned with the above two-term
summation formula \eqref{Slater} for ${}_2\psi_2$. Andrews
\cite{Andrews} established a three-term transformation formula which
is the key to proving many of Ramanujan's identities for partial
$\theta$-functions. In view of the symmetry in this formula, he
obtained a generalization of Ramanujan's $_1\psi_1$ summation:
\begin{eqnarray}\label{Andrews}
\lefteqn{d\sum_{n=0}^{\infty}\frac{(q/bc,acdf;q)_n}{(ad,df;q)_{n+1}}(bd)^n
-c\sum_{n=0}^{\infty}\frac{(q/bd,acdf;q)_n}{(ac,cf;q)_{n+1}}(bc)^n}\nonumber\\
&=&d\frac{(q,qd/c,c/d,abcd,acdf,bcdf;q)_\infty}{(ac,ad,bc,bd,cf,df;q)_\infty},
\quad\quad |bc|,|bd|<1.
\end{eqnarray}
Using the approach of parameter augmentation developed by Chen and
Liu \cite{Chen2}, we find that the two-term summation formula
\eqref{Slater} for ${}_2\psi_2$ series is a consequence of
 the
above identity \eqref{Andrews} of Andrews by bilateral extension and
parameter augmentation.

As is customary, we employ the notation and terminology of
basic hypergeometric series in
\cite{Gasper-Rahman}. For $|q|<1$, the $q$-shifted factorial is
defined by
$$(a;q)_\infty=
\prod_{k=0}^{\infty}(1-aq^k) \text{ and }(a;q)_n
=\frac{(a;q)_\infty}{(aq^n;q)_\infty}, \text{ for } n\in
\mathbb{Z}.$$ For convenience, we shall adopt the following notation
for multiple $q$-shifted factorials:
$$(a_1,a_2,\ldots,a_m;q)_n=(a_1;q)_n(a_2;q)_n\cdots(a_m;q)_n,$$
where $n$ is an integer or infinity. In particular, for a
nonnegative integer $k$, we have
\begin{equation} \label{x-k}
(a;q)_{-k}={1\over (aq^{-k};q)_{k}}.
\end{equation}

The (unilateral) basic hypergeometric series $_{r}\phi_s$ is defined
by
\begin{equation}
_{r}\phi_s\left[\begin{array}{cccccc}
a_1,&a_2,&\ldots,&a_r\\
b_1,&b_2,&\ldots,&b_s
\end{array};q,z\right]=\sum_{k=0}^{\infty}
\frac{(a_1,a_2,\ldots,a_r;q)_k}{(q,b_1,b_2,\ldots,b_s;q)_k}
\left[(-1)^kq^{k\choose2}\right]^{1+s-r}z^k,
\end{equation}
while the bilateral basic hypergeometric series $_{r}\psi_s$ is
defined by
\begin{equation}
_{r}\psi_s\left[\begin{array}{cccccc}
a_1,&a_2,&\ldots,&a_r\\
b_1,&b_2,&\ldots,&b_s
\end{array};q,z\right]=\sum_{k=-\infty}^{\infty}\frac{(a_1,a_2,\ldots,a_r;q)_k}
{(b_1,b_2,\ldots,b_s;q)_k}\left[(-1)^kq^{k\choose2}\right]^{s-r}z^k.
\end{equation}

\section{An Expansion Formula for the $_{2}\psi_2$ Series}

In this section, we  derive a representation for the $_{2}\psi_2$
series \eqref{2psi2-gene} in terms of two ${}_2\phi_1$ series. This
formula can be considered as a companion of Slater's formulas
\eqref{22-1} and \eqref{22-2}. We also present some consequences
including a two-term infinite product representation for the sum of
a well-poised $_2\psi_2$ series \eqref{2psi2-2}.

\begin{thm}\label{theorem2psi2-1}  We have
\begin{eqnarray}
&&{}_{2}\psi_2\left[\begin{array}{ccc}
a,&b\\
c,& d
\end{array};q,z\right]=\frac{(c/b,abz/d,dq/abz,q/d,q;q)_\infty}
{(c,az/d,q/a,q/b,cd/abz;q)_\infty}
{}_{2}\phi_{1}\left[\begin{array}{ccc}
cd/abz,&d/a\\
& dq/az
\end{array};q,\frac{bq}{d}\right]\nonumber\\
&&\quad\quad\quad\quad-\frac{(cq/d,b,d/a,az/q,q^2/az,q/d,q;q)_\infty}
{(d/q,c,bq/d,az/d,dq/az,q^2/d,q/a;q)_\infty}
{}_{2}\phi_{1}\left[\begin{array}{ccc}
aq/d,&bq/d\\
&cq/d
\end{array};q,z\right],\label{2psi2to2phi1}
\end{eqnarray}
where $|cd/ab|<|z|<1$ and $|bq/d|<1$.
\end{thm}
\begin{pf}
We start with a three-term transformation of ${}_{3}\phi_2$ series
\cite[Appendix III.33]{Gasper-Rahman}:
\begin{eqnarray*}
&&{}_{3}\phi_2\left[\begin{array}{ccc}
a, &b, &c\\
& d,& e
\end{array};q,\frac{de}{abc}\right]=
\frac{(e/b,e/c,cq/a,q/d;q)_\infty}
{(e,cq/d,q/a,e/bc;q)_\infty}
{}_{3}\phi_2\left[\begin{array}{ccc}
c, &d/a, &cq/e\\
& cq/a,& bcq/e
\end{array};q,\frac{bq}{d}\right]\nonumber\\
&&\quad-\frac{(q/d,eq/d,b,c,d/a,de/bcq,bcq^2/de;q)_\infty}
{(d/q,e,bq/d,cq/d,q/a,e/bc,bcq/e;q)_\infty}
{}_{3}\phi_2\left[\begin{array}{ccc}
aq/d, &bq/d, &cq/d\\
& q^2/d,& eq/d
\end{array};q,\frac{de}{abc}\right],
\end{eqnarray*}
where $|bq/d|,|de/abc|<1$.

Shifting the index of summation on the left hand side of the above
identity by $m$ such that the new sum runs from $-m$ to infinity,
and then replacing $a$, $b$, $d$, $e$ by $aq^{-m}$, $bq^{-m}$,
$dq^{-m}$, $eq^{-m}$, respectively, we get
\begin{eqnarray}
\lefteqn{\sum_{k=-m}^{\infty}\frac{(a,b,cq^m;q)_k}
{(q^{m+1},d,e;q)_k}\left(\frac{de}{abc}\right)^k
=\frac{(cq/e,q/d,q;q)_m}{(c,q/a,q/b;q)_m}
\frac{(e/b,e/c,cq^{1+m}/a,q^{1+m}/d;q)_\infty}
{(e,cq^{1+m}/d,q^{1+m}/a,e/bc;q)_\infty} }\nonumber\\
&&\times{}_{3}\phi_2\left[\begin{array}{ccc}
c, &d/a, &cq^{1+m}/e\\
& cq^{1+m}/a,& bcq/e
\end{array};q,\frac{bq}{d}\right]-\frac{(bcq^2/de,q/d,q;q)_m}
{(q^2/d,q/a,c;q)_m}\frac{(q^{1+m}/d,eq/d,b;q)_\infty}
{(d/q,e,bq/d;q)_\infty}\nonumber\\
&&\times\frac{(c,d/a,de/bcq,bcq^{2+m}/de;q)_\infty}
{(cq^{1+m}/d,q^{1+m}/a,e/bc,bcq/e;q)_\infty}
{}_{3}\phi_2\left[\begin{array}{ccc}
aq/d, &bq/d, &cq^{1+m}/d\\
& q^{2+m}/d,& eq/d
\end{array};q,\frac{de}{abc}\right],\label{3phi2-1}
\end{eqnarray}
where $|bq/d|,|de/abc|<1$.

Setting $m\rightarrow\infty$ in \eqref{3phi2-1} and assuming
$|c|<1$, Tannery's theorem \cite{Boas} enables us to interchange the
limit and the summation. This gives
\begin{eqnarray}
\lefteqn{{}_{2}\psi_2\left[\begin{array}{ccc}
a,&b\\
d,& e
\end{array};q,\frac{de}{abc}\right]
=\frac{(cq/e,q/d,q,e/b,e/c;q)_\infty} {(c,q/a,q/b,e,e/bc;q)_\infty}
{}_{2}\phi_1\left[\begin{array}{ccc}
c, &d/a\\
& bcq/e
\end{array};q,\frac{bq}{d}\right]}\nonumber\\
&&\quad-\frac{(bcq^2/de,q/d,q,eq/d,b,d/a,de/bcq;q)_\infty}
{(q^2/d,q/a,d/q,e,bq/d,e/bc,bcq/e;q)_\infty}
{}_{2}\phi_1\left[\begin{array}{ccc}
aq/d, &bq/d\\
& eq/d
\end{array};q,\frac{de}{abc}\right],\label{3phi2-2}
\end{eqnarray}
where $|bq/d|,|c|,|de/abc|<1$.

By the substitutions  $c \rightarrow de/abz$ and  $e \rightarrow c$
in \eqref{3phi2-2}, we get the desired formula. \qed
\end{pf}

Note that Theorem  \ref{theorem2psi2-1}  may be considered as a
bilateral extension of the following three-term transformation
formula \cite[Appendix III.31]{Gasper-Rahman}
\begin{eqnarray}\label{three-1}
&&{}_{2}\phi_1\left[\begin{array}{ccc}
a,&b\\
& d
\end{array};q,z\right]=\frac{(abz/d,q/d;q)_\infty}
{(az/d,q/a;q)_\infty} {}_{2}\phi_{1}\left[\begin{array}{ccc}
d/a,&dq/abz\\
& dq/az
\end{array};q,\frac{bq}{d}\right]\nonumber\\
&&\quad\quad\quad\quad\quad-\frac{(b,d/a,az/q,q^2/az,q/d;q)_\infty}
{(d/q,bq/d,az/d,dq/az,q/a;q)_\infty}
{}_{2}\phi_{1}\left[\begin{array}{ccc}
aq/d,&bq/d\\
&q^2/d
\end{array};q,z\right],
\end{eqnarray}
where $|bq/d|,|z|<1$. It is clear that  (\ref{three-1}) is a special
case of (\ref{2psi2to2phi1}) for  $c=q$.

Since Slater's formula \eqref{22-2} and our  formula
\eqref{2psi2to2phi1} deal with the same series, we are naturally led
to an identity on ${}_2\phi_1$ series. The right hand sides of
\eqref{22-2} and \eqref{2psi2to2phi1} give rise to the following
identity by replacing $a$, $b$, $c$, $z$ by $d/b$, $dz/q$, $adz/c$,
$bq/c$, respectively,
\begin{eqnarray}\label{2phi1w}
&&{}_{2}\phi_{1}\left[\begin{array}{ccc}
a,&b\\
&c
\end{array};q,z\right]=\frac{(abz/c,q/c;q)_\infty}
{(az/c,q/a;q)_\infty} {}_{2}\phi_{1}\left[\begin{array}{ccc}
cq/abz,&c/a\\
&cq/az
\end{array};q,\frac{bq}{c}\right]\nonumber\\[6pt]
&&\quad\quad\quad\quad\quad\quad\quad
+\left(\frac{q(1-a)(b,q/z,d/aq,aq^2/d,cq/adz,adz/c,q/c;q)_\infty}
{d(d,c/az,1/a,aq/c,dz/c,cq/dz,q/d;q)_\infty}\right.\nonumber\\[6pt]
&&+\left.\frac{(azq/c,dz/q,b,d/c,cq/d,q^2/dz,a;q)_\infty}
{(d/q,z,c,q^2/d,aq/c,dz/c,cq/dz;q)_\infty}\right)
{}_{2}\phi_{1}\left[\begin{array}{ccc}
q/b,&z\\
&azq/c
\end{array};q,\frac{bq}{c}\right].
\end{eqnarray}
It is worth noting that the parameter $d$  occurs only in the
factors of the second term on the right hand side of \eqref{2phi1w}.
Hence the sum of the two products in the parentheses does not depend
on $d$. This fact does not seem to be obvious by direct
verification. Setting $d=aq$, it follows that
\begin{eqnarray}\label{simi}
&&{}_{2}\phi_{1}\left[\begin{array}{ccc}
a,&b\\
&c
\end{array};q,z\right]=\frac{(abz/c,q/c;q)_\infty}
{(az/c,q/a;q)_\infty} {}_{2}\phi_{1}\left[\begin{array}{ccc}
cq/abz,&c/a\\
&cq/az
\end{array};q,\frac{bq}{c}\right]\nonumber\\[6pt]
&&\qquad +\frac{(az,b,c/a,q/az;q)_\infty} {(z,c,q/a,c/az;q)_\infty}
{}_{2}\phi_{1}\left[\begin{array}{ccc} q/b,&z\\
 &azq/c
\end{array};q,\frac{bq}{c}\right].
\end{eqnarray}
From  Heine's transformation \cite[Appendix III.1]{Gasper-Rahman}
\begin{equation}\label{fHeine}
{}_{2}\phi_1\left[\begin{array}{cc}
a,&b\\
&c
\end{array};q,z\right]=\frac{(b,az;q)_\infty}
{(c,z;q)_\infty}{}_{2}\phi_1\left[\begin{array}{cc}
c/b,&z\\
&az
\end{array};q,b\right],
\end{equation}
it is easily seen that \eqref{simi} is equivalent to \eqref{three-1} by the substitution $c \rightarrow d$.

\begin{cor}\label{theorem2psi2-2} We have
\begin{eqnarray}
&&{}_{2}\psi_2\left[\begin{array}{ccc}
a,&b\\
c,& d
\end{array};q,z\right]=\frac{(abz/d,c/b,dq/abz,q/d,q;q)_\infty}
{(c,az/d,q/a,q/b,cd/abz;q)_\infty}
{}_{2}\phi_{1}\left[\begin{array}{ccc}
cd/abz,&d/a\\
& dq/az
\end{array};q,\frac{bq}{d}\right]\nonumber\\
&&\quad\quad\quad\quad\quad+\frac{(d/a,b,az,q/az,q;q)_\infty}
{(d,c,d/az,z,q/a;q)_\infty} {}_{2}\phi_{1}\left[\begin{array}{ccc}
c/b,&z\\
&azq/d
\end{array};q,\frac{bq}{d}\right],\label{2psi2to2phi1cor}
\end{eqnarray}
where $|cd/ab|<|z|<1$ and $|bq/d|<1$.
\end{cor}
\begin{pf}
By  Heine's transformation \eqref{fHeine}, the second term on the right hand side of
\eqref{2psi2to2phi1} equals
\begin{eqnarray}
\lefteqn{-\frac{(cq/d,b,d/a,az/q,q^2/az,q/d,q;q)_\infty}
{(d/q,c,bq/d,az/d,dq/az,q^2/d,q/a;q)_\infty}
{}_{2}\phi_{1}\left[\begin{array}{ccc}
aq/d,&bq/d\\
&cq/d
\end{array};q,z\right]}\nonumber\\
&=&-\frac{(b,d/a,az/q,q^2/az,q/d,q,azq/d;q)_\infty}
{(d/q,c,az/d,dq/az,q^2/d,q/a,z;q)_\infty}
{}_{2}\phi_{1}\left[\begin{array}{ccc}
c/b,&z\\
&azq/d
\end{array};q,\frac{bq}{d}\right]\nonumber\\
&=&-\frac{(d/a,b,az,q/az,q;q)_\infty} {(d,c,d/az,z,q/a;q)_\infty}
\frac{(1-az/q)(1-q/d)(1-d/az)}{(1-d/q)(1-az/d)(1-q/az)}\nonumber\\
&& \times {}_{2}\phi_{1}\left[\begin{array}{ccc}
c/b,&z\\
&azq/d
\end{array};q,\frac{bq}{d}\right]\nonumber\\
&=&\frac{(d/a,b,az,q/az,q;q)_\infty} {(d,c,d/az,z,q/a;q)_\infty}
{}_{2}\phi_{1}\left[\begin{array}{ccc}
c/b,&z\\
&azq/d
\end{array};q,\frac{bq}{d}\right].
\end{eqnarray}\qed
\end{pf}
\begin{rem}
 Corollary \ref{theorem2psi2-2} can also be obtained from
the following three-term transformation formula \cite[Appendix
III.34]{Gasper-Rahman}
\begin{eqnarray*}
&&{}_{3}\phi_2\left[\begin{array}{ccc}
a, &b, &c\\
& d,& e
\end{array};q,\frac{de}{abc}\right]=
\frac{(e/b,e/c;q)_\infty} {(e,e/bc;q)_\infty}
{}_{3}\phi_2\left[\begin{array}{ccc}
d/a, &b, &c\\
& d,& bcq/e
\end{array};q,q\right]\nonumber\\
&&\quad+\frac{(d/a,b,c,de/bc;q)_\infty} {(d,e,bc/e,de/abc;q)_\infty}
{}_{3}\phi_2\left[\begin{array}{ccc}
e/b, &e/c, &de/abc\\
&de/bc,&eq/bc
\end{array};q,q\right].
\end{eqnarray*}
Shifting the summation index by $m$ on the left hand side and
replacing $a$, $c$, $d$, $e$ by $aq^{-m}$, $cq^{-m}$, $dq^{-m}$,
$eq^{-m}$, respectively, we are led to \eqref{2psi2to2phi1cor} by
taking the limit  $m\rightarrow\infty$ and making suitable
substitutions.
\end{rem}

As a consequence of Corollary \ref{theorem2psi2-2}, we may deduce
the following
 expansion of a
$_2\psi_2$ series in terms of a $_2\phi_1$ series \cite[Eq.
(3.13.1.7)]{Exton}. Setting $z=q/a$ in \eqref{2psi2to2phi1cor}, the
second summation on the right hand side vanishes. It follows from
\eqref{fHeine} that
\begin{equation}
{}_{2}\psi_2\left[\begin{array}{cc}
a,&b\\
c,&d
\end{array};q,\frac{q}{a}\right]=
\frac{(c/b,d/b,bq/a,q;q)_{\infty}}{(c,d,q/a,q/b;q)_{\infty}}
{}_{2}\phi_1\left[\begin{array}{cc}
bq/c,&bq/d\\
&bq/a
\end{array};q,\frac{cd}{bq}\right],
\end{equation}
which was originally derived from a double sum transformation
formula of Slater, see \cite[Section 3.13]{Exton}.

\begin{cor}\label{2psi2wellpoised} We have
\begin{eqnarray}
&&{}_{2}\psi_{2}\left[\begin{array}{ccc}
b,&c\\
aq/b,&aq/c
\end{array};q,-\frac{aq}{bc}\right]=
\frac{(-b,aq/bc,-q/b,b/a,q;q)_\infty(aq^2/c^2;q^2)_\infty}
{(aq/c,-1,q/c,q/b,-aq/bc;q)_\infty(b^2/a;q^2)_\infty}\nonumber\\
&&\quad\quad\quad\quad\quad\quad\quad\quad
+\frac{(aq/bc,b,-aq/b,-b/a,q;q)_\infty(aq^2/c^2;q^2)_\infty}
{(aq/b,aq/c,-1,-aq/bc,q/c;q)_\infty(b^2/a;q^2)_\infty}
,\label{2psi22term}
\end{eqnarray}
where $|aq/bc|<1$.
\end{cor}
\begin{pf}
Setting $c=cq/a$, $d=cq/b$, and $z=-cq/ab$ in
\eqref{2psi2to2phi1cor}, we find that the summations on the right
hand side of the identity are both equal to
\begin{equation}
\sum_{k=0}^{\infty}\frac{(c^2q^2/a^2b^2;q^2)_k}{(q^2;q^2)_k}\left(\frac{b^2}{c}\right)^k,
\end{equation}
which can be summed by the Cauchy $q$-binomial theorem
\cite[Appendix II.3]{Gasper-Rahman}
\begin{equation}\label{Cauchy}
\sum_{n=0}^{\infty}\frac{(a;q)_n}{(q;q)_n}z^n
=\frac{(az;q)_\infty}{(z;q)_\infty},\quad\quad\quad\quad |z|<1.
\end{equation}
Thus the following relation holds
\begin{eqnarray*}
{}_{2}\psi_{2}\left[\begin{array}{ccc}
a,&b\\
cq/a,&cq/b
\end{array};q,-\frac{cq}{ab}\right]&=&
\frac{(-b,cq/ab,-q/b,b/c,q;q)_\infty(cq^2/a^2;q^2)_\infty}
{(cq/a,-1,q/a,q/b,-cq/ab;q)_\infty(b^2/c;q^2)_\infty}\nonumber\\
&&+\frac{(cq/ab,b,-cq/b,-b/c,q;q)_\infty(cq^2/a^2;q^2)_\infty}
{(cq/b,cq/a,-1,-cq/ab,q/a;q)_\infty(b^2/c;q^2)_\infty} .
\end{eqnarray*}
The proof is thus completed by interchanging $a$ and $c$.\qed
\end{pf}

Combining \eqref{2psi22term} and \eqref{2psi21term}, we are led to
the following identity
\begin{eqnarray}
\lefteqn{(-b,-q/b,b/a,aq/b;q)_\infty+(b,q/b,-b/a,-aq/b;q)_\infty}\qquad
\nonumber\\[6pt]
&=&\frac{2(aq,q/a,b^2/a,aq^2/b^2;q^2)_\infty}{(q;q^2)_\infty^2}.\label{infprod}
\end{eqnarray}
To restate the above identity in a symmetric form, we replace $a$ by
$b/a$ in \eqref{infprod}.

\begin{thm}\label{infprodthm} We have
\begin{eqnarray}
(a,-b,q/a,-q/b;q)_\infty+(-a,b,-q/a,q/b;q)_\infty
=\frac{2(ab,q^2/ab,aq/b,bq/a;q^2)_\infty}{(q;q^2)_\infty^2}.
\end{eqnarray}
\end{thm}

More identities on sums of infinite products have been found by
Bailey \cite{Bailey1} and Slater \cite{slater1,slater2,slater3}.

While no attempt will be made to derive a closed product formula for
the series (\ref{2psi2gauss}), we obtain a formula involving a
product and a summation which has the advantage that it reduces to
the $q$-Gauss summation \eqref{Gauss} when $c=q$ or $d=q$. Combining
Corollary \ref{theorem2psi2-2} and  Cauchy's $q$-binomial theorem
\eqref{Cauchy}, we deduce

\begin{cor}
\begin{eqnarray}
&&{}_{2}\psi_2\left[\begin{array}{ccc}
a,&b\\
c,& d
\end{array};q,\frac{cd}{abq}\right]=\frac{(c/b,c/q,q^2/c,q/d;q)_\infty}
{(c,c/bq,q/a,q/b;q)_\infty}\sum_{k=0}^{\infty}\frac{(d/a;q)_k}{(bq^2/c;q)_k}
\left(\frac{bq}{d}\right)^k\nonumber\\
&&\quad\quad\quad\quad\quad\quad\quad\quad\quad\quad+\frac{(c/a,d/a,b,cd/bq,bq^2/cd,q;q)_\infty}
{(c,d,bq/c,bq/d,q/a,cd/abq;q)_\infty},\label{biGauss}
\end{eqnarray}
where $|bq/d|,|cd/abq|<1$.
\end{cor}

\section{A Two-term Summation Formula for  $_{2}\psi_2$}

In this section, we show that a two-term summation formula for the
${}_2\psi_2$ series \eqref{Slater} due to Slater can be derived from
an identity of Andrews \eqref{Andrews} by bilateral extension and parameter augmentation.

We recall that the $q$-difference operator, or Euler derivative, is
defined as
\begin{equation}
D_q\{f(a)\}=\frac{f(a)-f(aq)}{a}.
\end{equation}
The $q$-shift operator $\eta$ in the literature
\cite{Andrews1,Roman} is defined as follows:
\begin{equation}
\eta\{f(a)\}=f(aq) \hskip 5mm \text{and} \hskip 5mm
\eta^{-1}\{f(a)\}=f(aq^{-1}),
\end{equation}
which was introduced by Rogers in \cite{Rogers1,Rogers2,Rogers3}.

In \cite{Roman}, Roman combined $q$-differential operator and the
$q$-shift operator to built an operator which was denoted by
$\theta$ in \cite{Chen2}:
\begin{equation}
\theta=\eta^{-1}D_q.
\end{equation}
In \cite{Chen2}, Chen and Liu introduced the operator:
\begin{equation}\label{Eth}
E(b\theta)=\sum_{n=0}^{\infty}\frac{(b\theta)^n q^{{n \choose
2}}}{(q;q)_n},
\end{equation}
and proved the following basic relations:
\begin{eqnarray}
E(b\theta)\left\{(at;q)_\infty\right\}&=&(at,bt;q)_\infty,\label{o1}\\[6pt]
E(b\theta)\left\{(as,at;q)_\infty\right\}&=&
\frac{(as,at,bs,bt;q)_\infty}{(abst/q;q)_\infty},\quad\quad
|abst/q|<1.\label{o2}
\end{eqnarray}
The procedure to apply the operator $E(b\theta)$ in order to derive
a new identity is called parameter augmentation.

The following theorem is equivalent to Slater's formula
\eqref{Slater}, as pointed out by Ismail and Rahman
\cite{Ismail-Rahman}. We proceed to demonstrate how to derive it
from the identity (\ref{Andrews}) of Andrews by bilateral extension
and parameter augmentation.

\begin{thm}\label{2psi2} We have
\begin{eqnarray}\label{2psi2g}
&&{}_{2}\psi_2\left[\begin{array}{cc}
a,&b\\
c,& d
\end{array};q,\frac{cd}{abq}\right]-\frac{\alpha}{q}\frac{(q/c,q/d,\alpha/a,\alpha/b;q)_\infty}
{(q/a,q/b,\alpha/c,\alpha/d;q)_\infty}
{}_{2}\psi_2\left[\begin{array}{cc}
aq/\alpha,&bq/\alpha\\
cq/\alpha,& dq/\alpha
\end{array};q,\frac{cd}{abq}\right]\nonumber\\[5pt]
&&\quad\quad\quad=\frac{(\alpha,q/\alpha,cd/\alpha q, \alpha
q^2/cd,q,c/a,c/b,d/a,d/b;q)_\infty} {(c/\alpha,\alpha
q/c,d/\alpha,\alpha q/d,c,d,q/a,q/b,cd/abq;q)_\infty},
\end{eqnarray}
where $|cd/abq|<1$.
\end{thm}

\begin{pf}
Shifting the index of summation by $m$ and then replacing $a$, $b$,
$f$ by $aq^{-m}$, $bq^{m}$, $fq^{-m}$ in \eqref{Andrews},
respectively, we obtain
\begin{eqnarray}
\lefteqn{\frac{d(q^{1-m}/bc,acdfq^{-2m};q)_m\left(bdq^m\right)^m}
{(1-adq^{-m})(1-dfq^{-m})(adq^{1-m},dfq^{1-m};q)_m}
\sum_{k=-m}^{\infty}
\frac{(q/bc,acdfq^{-m};q)_k}{(adq,dfq;q)_k}\left(bdq^m\right)^k}\nonumber\\
&&-\frac{c(q^{1-m}/bd,acdfq^{-2m};q)_m\left(bcq^m\right)^m}
{(1-acq^{-m})(1-cfq^{-m})(acq^{1-m},cfq^{1-m};q)_m}
\sum_{k=-m}^{\infty}
\frac{(q/bd,acdfq^{-m};q)_k}{(acq,cfq;q)_k}\left(bcq^m\right)^k\nonumber\\
&=&\frac{d(q,qd/c,c/d,abcd,acdfq^{-2m},bcdf;q)_\infty}
{(acq^{-m},adq^{-m},bcq^{m},bdq^{m},cfq^{-m},dfq^{-m};q)_\infty}.\label{Andrewm}
\end{eqnarray}
Letting $m\rightarrow\infty$ in \eqref{Andrewm} and employing
Tannery's theorem, we get
\begin{eqnarray}
\lefteqn{\frac{c(bc;q)_\infty}{(1/ad,1/df;q)_\infty}\sum_{k=-\infty}^{\infty}
\frac{(q/bc;q)_k}{(adq,dfq;q)_k}\left(-abcd^2f\right)^kq^{k\choose2}}\nonumber\\
&&\quad\quad-\frac{d(bd;q)_\infty}{(1/ac,1/cf;q)_\infty}\sum_{k=-\infty}^{\infty}
\frac{(q/bd;q)_k}{(acq,cfq;q)_k}\left(-abc^2df\right)^kq^{k\choose2}\nonumber\\
&=&\frac{acd^2f(q,qd/c,c/d,abcd,acdf,bcdf,q/acdf;q)_\infty}
{(ac,ad,cf,df,q/ac,q/ad,q/cf,q/df;q)_\infty}.\label{Andrewminf}
\end{eqnarray}
Now, \eqref{Andrewminf} can be written as
\begin{eqnarray}
\lefteqn{\frac{c}{(1/ad,1/df;q)_\infty}\sum_{k=-\infty}^{\infty}
\frac{(bcq^{-k};q)_\infty}{(adq,dfq;q)_k}\left(ad^2fq\right)^kq^{2{k\choose2}}}\nonumber\\
&&\quad-\frac{d}{(1/ac,1/cf;q)_\infty}\sum_{k=-\infty}^{\infty}
\frac{(bdq^{-k};q)_\infty}{(acq,cfq;q)_k}\left(ac^2fq\right)^kq^{2{k\choose2}}\nonumber\\
&=&\frac{acd^2f(q,qd/c,c/d,abcd,acdf,bcdf,q/acdf;q)_\infty}
{(ac,ad,cf,df,q/ac,q/ad,q/cf,q/df;q)_\infty}.\label{Andrewminfd}
\end{eqnarray}
Next, applying $E(g\theta)$ to both sides of \eqref{Andrewminfd}
with respect to the parameter $b$ gives
\begin{eqnarray}
\lefteqn{\frac{c}{(1/ad,1/df;q)_\infty}\sum_{k=-\infty}^{\infty}
\frac{\left(ad^2fq\right)^kq^{2{k\choose2}}}{(adq,dfq;q)_k}
E(g\theta)\left\{(bcq^{-k};q)_\infty\right\}}\nonumber\\
&&\quad\quad-\frac{d}{(1/ac,1/cf;q)_\infty}\sum_{k=-\infty}^{\infty}
\frac{\left(ac^2fq\right)^kq^{2{k\choose2}}}{(acq,cfq;q)_k}
E(g\theta)\left\{(bdq^{-k};q)_\infty\right\}\nonumber\\
&=&\frac{acd^2f(q,qd/c,c/d,acdf,q/acdf;q)_\infty}
{(ac,ad,cf,df,q/ac,q/ad,q/cf,q/df;q)_\infty}
E(g\theta)\left\{abcd,bcdf;q)_\infty\right\}.\label{Andrewoperator}
\end{eqnarray}
From \eqref{o1} and \eqref{o2}, it is evident that
\begin{equation}
E(g\theta)\left\{(bcq^{-k};q)_\infty\right\}=(bcq^{-k},cgq^{-k};q)_\infty,\label{action1}
\end{equation}
\begin{equation}
E(g\theta)\left\{(bdq^{-k};q)_\infty\right\}=(bdq^{-k},dgq^{-k};q)_\infty,\label{action2}
\end{equation}
and
\begin{equation}
E(g\theta)\left\{abcd,bcdf;q)_\infty\right\}=
\frac{(abcd,acdg,bcdf,cdfg;q)_\infty}{(abc^2d^2fg/q;q)_\infty}\label{action3}.
\end{equation}
Substituting \eqref{action1}, \eqref{action2}, and \eqref{action3}
into \eqref{Andrewoperator}, we see that
\begin{eqnarray}
\lefteqn{\frac{c(bc,cg;q)_\infty}{(1/ad,1/df;q)_\infty}
{}_{2}\psi_2\left[\begin{array}{ccc}
q/bc,&q/cg\\
adq,& dfq
\end{array};q,\frac{abc^2d^2fg}{q}\right]}\nonumber\\
&&\quad\quad-\frac{d(bd,dg;q)_\infty}{(1/ac,1/cf;q)_\infty}
{}_{2}\psi_2\left[\begin{array}{ccc}
q/bd,&q/dg\\
acq,& cfq
\end{array};q,\frac{abc^2d^2fg}{q}\right]\nonumber\\
&=&\frac{acd^2f(q,qd/c,c/d,abcd,acdf,acdg,bcdf,q/acdf,cdfg;q)_\infty}
{(ac,ad,cf,df,q/ac,q/ad,q/cf,q/df,abc^2d^2fg/q;q)_\infty},\label{2psi22termpsi}
\end{eqnarray}
where $|abc^2d^2fg/q|<1$.

Finally, the proof is completed by replacing $a$, $b$, $c$, $d$,
$f$, $g$ by $c/fq$, $e$, $q/ae$, $f$, $d/fq$, $ae/b$, respectively,
and then setting $aef=\alpha$. \qed
\end{pf}

Substitute $a$, $b$, $c$, $d$, $\alpha$ with $qa/e$, $qb/e$, $qc/e$,
$q^2/e$, $fq/e$ in \eqref{2psi2g}, respectively, we may recover the
original formula \eqref{Slater} due to Slater.

If we set $d=q$ in \eqref{2psi2g}, then the second term on the left
hand side
 vanishes, and so we get the
$q$-Gauss summation \eqref{Gauss} as a special of \eqref{2psi2g}.

To conclude this paper, we represent (\ref{2psi2g}) in an equivalent form
and give the explicit substitutions to reach Askey's $q$-extension of Cauchy's
beta integral  \eqref{integral}.
By the relation
\begin{eqnarray}
{}_{2}\psi_2\left[\begin{array}{ccc}
a,&b\\
c,&d
\end{array};q,z\right]={}_{2}\psi_2\left[\begin{array}{ccc}
q/c,&q/d\\
q/a,&q/b
\end{array};q,\frac{cd}{abz}\right],
\end{eqnarray}
we may rewrite \eqref{2psi2g} as
\begin{eqnarray}\label{22}
&&\frac{(q/a,q/b;q)_\infty}{(q/c,q/d;q)_\infty}{}_{2}\psi_2\left[\begin{array}{cc}
q/c,&q/d\\
q/a,&q/b
\end{array};q,q\right]-\frac{\alpha}{q}\frac{(\alpha/a,\alpha/b;q)_\infty}
{(\alpha/c,\alpha/d;q)_\infty}{}_{2}\psi_2\left[\begin{array}{cc}
\alpha/c,&\alpha/d\\
\alpha/a,&\alpha/b
\end{array};q,q\right]\nonumber\\
&&\quad\quad\quad\quad =\frac{(\alpha,q/\alpha,cd/\alpha q, \alpha
q^2/cd,q,c/a,c/b,d/a,d/b;q)_\infty} {(c/\alpha,\alpha
q/c,d/\alpha,\alpha q/d,c,d,q/c,q/d,cd/abq;q)_\infty},
\end{eqnarray}
where $|cd/abq|<1$. Replacing $a$, $b$, $c$, $d$, $\alpha$ by $q/c$,
$-q/d$, $q/a$, $-q/b$, $q$, respectively, then \eqref{22} takes the
form of Askey's $q$-extension of Cauchy's beta integral.

%--------------------------------------------------------------------
\vspace{.2cm} \noindent{\bf Acknowledgments.} We are grateful to George Andrews and Richard Askey for their valuable comments. This work was
supported by the 973 Project on Mathematical Mechanization, the
Ministry of Science and Technology, the Ministry of Education and
the National Science Foundation of China.
%-------------

\end{document}